\documentclass[12pt]{article}
\usepackage{amssymb}
\newtheorem{theorem}{\sc Theorem}[section]
\newtheorem{definition}{\sc Definition}[section]
\newtheorem{lemma}{\sc Lemma}[section]

\newtheorem{example}{\sc Example}[section]
\newtheorem{remark}{\it Remark}[section]
\newtheorem{Ack}{\it Acknowledgements}[section]
\begin{document}
\baselineskip=24 pt
\title{\bf On the  reduction principle for the  hybrid equation}
\author{ M. U. Akhmet\thanks{M.U. Akhmet is previously known as M. U. Akhmetov. }}

\date{{\small Department of Mathematics and Institute of Applied Mathematics, Middle East
Technical University, 06531 Ankara, Turkey}}
\maketitle

\noindent {\bf Address:} M. Akhmet,
Department of Mathematics, Middle East
Technical University, 06531 Ankara, Turkey,\\
{\bf fax:} 90-312-210-12-82\\
\noindent {\bf e-mail:} marat@metu.edu.tr
 
\vspace{0.5cm}

\newpage
\noindent {\bf Proposed running head:} The Reduction Principle
\begin{abstract}
In this  paper we introduce a new type of differential equations with piecewise constant argument (EPCAG), more general than EPCA  \cite{cw1, w}. The Reduction Principle \cite{pliss} is proved for  EPCAG. The continuation of solutions is investigated.
We establish the existence of global integral manifolds of quasilinear EPCAG and investigate the stability of the zero solution.  Since  the method of reduction to discrete equations \cite{cw1} is difficult to apply to EPCAG, a new technique of investigation of equations with piecewise argument, based on an integral representation formula, is proposed. The approach  can be fruitfully applied for investigation  stability, oscillations, controllability and many other problems of EPCAG.

\end{abstract}
\maketitle
\section{Introduction and Preliminaries}
\subsection{Definitions and the description of the system} 
The theory of integral manifolds was founded by  H. Poincar\'e and A. M. Lyapunov \cite{poin, lya}, and it became a very powerful instrument for investigating various problems of the qualitative theory of differential equations. For the last several decades, many researchers have been studying the methods of reducing high  dimensional problems to low dimensional ones. When discussing this problem for long-time dynamics of differential equations, we should consider  the Reduction Principle  \cite{pliss, pliss1}. One can read about the history of the principle   in \cite {l, malkin,pliss} and papers cited there. The principle  was utilized  in the center manifold theory, as well as in the theory of inertial manifolds \cite{carr, henry,foias}. It is natural that the exploration  of the properties and neighborhoods of manifolds is one of the most interesting problems of the theory of differential equations \cite{b1,carr,chow,h,h1,kelley, lat, pal1,pugh}. One should  not be surprised that  manifolds and the reduction principle are  one of the major subjects of investigation for specific types of differential  and difference equations \cite{akhmet, aul, chic,ch,diek, h,h1,hart,henry,foias,lee,p1,pugh,sp,stokes}. Our main goal in this paper is to extend the principle to  the differential equations with piecewise constant argument
of generalized type. For this purpose, we have developed another approach to the investigation, different from what was proposed by  the founders of the EPCA theory \cite{cw1,w}. 

Let $\mathbb Z, \mathbb N$ and $\mathbb R$ be the sets of all integers, natural and real
numbers, respectively. Denote by $||\cdot||$ the Euclidean norm in
$\mathbb R^n$, $n \in \mathbb N.$
Fix two real-valued sequences $\theta_i, \zeta_i, i \in \mathbb Z,$ such that $   \theta_i < \theta_{i+1}, \theta_i\leq \zeta_i \leq \theta_{i+1}$ for all $i \in \mathbb Z,$ 
$|\theta_i| \rightarrow \infty$ as $ |i| \rightarrow \infty,$ and there exists a number $\theta >0$ such that $ \theta_{i+1} - \theta_{i} \leq \theta, i \in \mathbb Z.$  
In this paper we are concerned with the quasilinear system
\begin{eqnarray}
&& z' = Az + f(t,z(t),z(\beta(t))),
\label{1}
\end{eqnarray}
where $z \in \mathbb R^n, t \in \mathbb R, \beta(t) =\zeta_i,$ 
if $t \in [\theta_i, \theta_{i+1}), i \in \mathbb Z.$ 
One can easily see that equation (\ref{1}) has the form 
\begin{eqnarray}
	&& z' = Az + f(t,z(t), \bar z),
	\label{1'}
\end{eqnarray}
 if $t \in [\theta_i, \theta_{i+1}),  \bar z = z(\zeta_i), i \in \mathbb Z.$ 

The theory of differential equations with  piecewise constant argument\\ (EPCA) of the type
\begin{eqnarray}
&& \frac{dx(t)}{dt} = f(t,x(t),x([ t])),
\label{e2}
\end {eqnarray}
where $[\cdot]$ signifies the greatest integer function,  was initiated in  \cite{cw1} and 
has been developed by many authors
 \cite{aa1, cwm, Gy, ky,p1,w,wl}. They are hybrid equations, in that they combine
 the properties of both continuous systems and  discrete equations. For example, even a scalar logistic equation may produce chaos \cite{Hop,g,K}, when the solutions are continuous functions.
 
The novel idea of our paper is that   system (\ref{1})  is EPCA of general type
 (EPCAG)  for equation (\ref{e2}). 
Indeed if we take $\zeta_i=\theta_i = i, i \in \mathbb Z,$ then   (\ref{1}) takes the form of (\ref{e2}).
Another EPCA which can be easily written as EPCAG is the equation 
alternately of retarded and advanced type \cite{cwm,w}
\begin{eqnarray}
&& \frac{dx(t)}{dt} = f(t,x(t),x(2[\,(t+1)/2\,])).
\label{e29}
\end {eqnarray}
One can check that (\ref{1}) takes  the form of (\ref{e29}) if 
$\theta_i = 2i-1, \zeta_i = 2i, i \in \mathbb Z.$ Moreover, the system considered in 
\cite{akhmet2} is a particular case of (\ref{1}), too.

The following assumptions will be needed throughout the paper:
\begin{itemize}
\item[(C1)]  $A$ is a constant $n\times n$ real valued matrix;
\item[(C2)]  $f(t,x,z)$ is continuous in the first argument,  $f(t, 0,0) = 0, t \in \mathbb R,$ and $f$ is Lipschitzian in 
the second and  third arguments with a  positive Lipschitz
constant  $l$ such that
$$||f(t,z_1,w_1) -  f(t,z_2,w_2)|| \leq l ( ||z_1-z_2|| + ||w_1-w_2||)$$
for all $t \in \mathbb R$ and $z_1, z_2, w_1, w_2 \in \mathbb R^n.$
 
\item[(C3)]  If we denote by $\lambda_j, j= \overline{1,n},$ the eigenvalues of  matrix $A,$ then there exists a positive integer $k$ such that $\mu = \max_{j=\overline{1,k}} \Re \lambda_j <0,$ 
and $\min_{j=\overline{k+1,n}} \Re \lambda_j =0,$
where $\Re \lambda_j$ denotes the real part of the eigenvalue $\lambda_j$ of  matrix $A.$

The previous condition implies that, without loss of generality, we can assume that
 \item[(C4)]  
\begin{displaymath}
A = \left( \begin{array}{cc}
B_+ & 0 \\
0 & B_-\end{array} \right),
\end{displaymath}
where square matrices $B_+$ and $B_-$ are of dimension $k$ and $n-k$ respectively, $\lambda_j, j =\overline{1,k},$
are eigenvalues of the matrix $B_-$ and   $\lambda_j, j =\overline{k+1,n},$ are the eigenvalues of  matrix $B_+.$
\end{itemize}
 
The existing method of investigation 
 of EPCA, as   proposed 
by founders, is based on  the reduction of EPCA to  discrete equations. 
It is obvious that this method is not applicable to the present problem. A new  approach is based on  the construction of an equivalent integral 
equation. Consequently, we prove  a corresponding equivalence lemma for every result of our paper. Thus, when investigating EPCAG, we need not impose any conditions on  the reduced discrete equations, and  hence require more easily verifiable conditions, similar to those for ordinary differential equations. It becomes  less cumbersome  to solve the problems  of EPCAG theory (as well as of   EPCA theory). 

The theory of 
EPCAG (EPCA) necessitates a  more careful discussion of the continuation problem. The subject of backward continuation for functional differential equations was considered in \cite{h1}.
In our paper it is necessary to analyze  the forward continuation, too, as we also deal with equations  alternately of retarded and advanced type. The backward continuation of the solutions of EPCA was investigated in \cite{cw1} through the solvability of  certain difference equations. For our needs,  we shall introduce less formal definitions than those in  \cite{cw1}, since we  consider  integral manifolds,  and  it is natural  to discuss the global  continuation of a solution of  (\ref{1}) as well as its uniqueness on these manifolds.  
\begin{definition} A solution  $z(t,t_0,z_0), \zeta_i < t_0 \leq \theta_{i+1},$
of (\ref{1})  is said to be backward continued to $t=\zeta_i$
if there exists a solution $z(t,\zeta_i,\bar z)$  of (\ref{1'})  
such that $z(t_0,\zeta_i,\bar z) =  z_0.$
The solution  $z(t,t_0,z_0)$ is uniquely backward
continued to  $t=\zeta_i$  if the continuation is unique.
\label{defn2}
\end{definition}
\begin{definition} A solution  $z(t,t_0,z_0), \theta_i \leq t_0 < \zeta_{i},$
of (\ref{1})  is said to be forward  continued to $t=\zeta_{i},$
if there exists a solution $z(t,\zeta_{i},\bar z)$  of (\ref{1'})  
such that $z(t_0,\zeta_{i},\bar z) =  z_0.$   The solution $z(t,t_0,z_0)$   is
uniquely continued to  $t=\zeta_{i}$  if the continuation is unique.
\label{defn+2}
\end{definition}
The following  example  shows that even for simple EPCAG the  continuation
of some solutions can  fail. 
\begin{example} 
Consider the following EPCA	
\begin{eqnarray}
&& z'= 3z - z^2(2[\,(t+1)/2\,]),
\label{ex1} 	
\end{eqnarray}
where $z \in \mathbb R,  t \in \mathbb R.$ 
Let us show that not all solutions of (\ref{ex1}) can be forward continued.
Consider the interval $[-1,0].$  Fix $z_0 \in \mathbb R,$ and let
 $z(t,0,z_0)$ be a solution of  (\ref{ex1}). It is clear that for a solution
  $z(t,-1,x_0)$ to be forward continued to $t=0,$ in the sense of Definition \ref{defn2}, the equation $[{\rm e}^3-1]z_0^2 -{\rm e}^3z_0 -x_0 = 0$ must be solvable with respect to $z_0.$  Since this equation can not be solved for all $z_0 \in \mathbb R,$ the assertion is proved. 
Further, we shall consider the uniqueness of the continuation.  Now let us focus  on the backward continuation.
Consider the interval $[0,1].$
Fix numbers $z_0, z_1 \in \mathbb R$ such that 
$(z_0 + z_1)(1-{\rm e}^3) = {\rm e}^3.$ 
Denote $z_0(t) = z(t,0,z_0)$ and $ z_1(t) = z(t,0,z_1)),$ which are  solutions of (\ref{ex1}).
 They are  the solutions of the equations 
$z'= 3z - z_0^2$ and $z'= 3z - z_1^2,$ respectively, and they are defined on $[0,1).$ 
Since 
$z_j(t) = {\rm e}^{3t}z_j + \int_0^t{\rm e}^{3(t-s)}z_j^2ds, j=0,1,$ one can obtain, using the continuity of solutions,  that
$z_0(1) =z_1(1).$ That is,  the solution $z(t,1,z_1(1))$ of (\ref{ex1}) can not be continued back to $t=0$ uniquely. 
\end{example}
Next we  consider the construction procedure  for a solution of an initial value problem.
We define the solution  only for decreasing $t,$ but one can easily see that the definition is similar for increasing $t.$ 

Let us assume  that $\theta_i\leq \zeta_i <t_0 \leq    \theta_{i+1}$  for some $i \in \mathbb Z.$
Suppose that
 $z(t)= z(t,t_0,z_0)$ is back continued from $t_0$ to $t= \zeta_i$ in the sense of Definition 
\ref{defn2}. Then conditions $(C1)$ and $(C2)$ imply that $z(t)$ can be continued to $t=\theta_i,$
as it is a solution of the following system of ordinary differential equations
$z' = Az + f(t,z(t),z(\zeta_i))$
on $[\theta_i, \theta_{i+1}).$

Next, we suppose that  $z(t,\theta_i,z(\theta_i,t_0,z_0))$ is back continued  from $t=\theta_i$ to 
$t= \zeta_{i-1}$ in the sense of Definition 
\ref{defn2}. Then again we can conclude that $z(t)$ can be continued to $t=\theta_{i-1}.$
If  $z(t,\theta_{i-1},z(\theta_{i-1},t_0,z_0))$ is back continued from $\theta_{i-1}$ to 
$t= \zeta_{i-2}$ in the sense of Definition 
\ref{defn2}, then $z(t)$ can be continued to $t=\theta_{i-2}.$
Proceeding in this way and assuming that $z(t,\theta_{j},z(\theta_{j},t_0,z_0))$ is back continued from $\theta_{j}$  to 
$t= \zeta_{j-1}$ in the sense of Definition 
\ref{defn2} for all $j \le i,$ we can find that 
$z(t, t_0, z_0)$ is back continuable to $- \infty.$
{\it We call this continuation  of   $z(t, t_0, z_0)$ to $- \infty$ as a solution of 
(\ref{1}) on $(-\infty, t_0].$} Similarly, one can define a  solution  as the continuation on
an interval $[t_0, \infty)$ as  time is increasing.
On the basis of the above discussion we can conclude that the following theorem is valid.
\begin{theorem} Assume that conditions $(C1)$ and $(C2)$ hold, $\zeta_i$ is the maximal among $\zeta_j$ which are 
smaller than $t_0.$  Then  $z(t, t_0, z_0)$ exists on $(- \infty,t_0]$ if and only if 
$z(t, \theta_j, z(\theta_j, t_0, z_0))$ is back continuable to $\zeta_{j-1}$ in the sense of Definition \ref{defn2} for all $ j \le i.$
\label{gu-gu}
\end{theorem}
A similar theorem can be proved for forward continuation.\\
In what follows we shall say  that a solution $z(t)$ {\it is continued} if it {\it is continued 
backward} or/and  {\it  forward}.
\begin{definition} A function $z(t)=z(t,t_0,z_0), z(t_0) = z_0,  \theta_i \le t_0 < \theta_{i+1}, i \in \mathbb Z,$
 is a solution of (\ref{1})
on the interval $[\theta_i, \infty)$  if
the following conditions are fulfilled:
\begin{enumerate}
\item[(i)]   z(t)  is continuable to $t=\zeta_i;$
\item[(ii)] the derivative $z'(t)$ exists at each point $t \in [\theta_i, \infty)$ with the possible exception
of the points $\,\theta_j \in [\theta_i, \infty)$   where one-sided derivatives exist;
\item[(iii)] equation (\ref{1}) is satisfied for $z(t)$ at each point $t \in [\theta_i, \infty)\backslash \{\theta_j\},$ and it holds for the right derivative of $z(t)$ at the points $\theta_j \in [\theta_i, \infty), \, j \in \mathbb Z.$
\end{enumerate}
\label{defn3}
\end{definition}
\begin{remark} One can see that Definition \ref{defn3} is  a slightly changed version of a definition from \cite{cw1}, adapted for our general case.
\end{remark}
\begin{definition} A function $z(t)=z(t,t_0,z_0), z(t_0) = z_0,  \theta_i < t_0 \leq \theta_{i+1}, i \in \mathbb Z,$
 is a solution of (\ref{1})
on the interval $(- \infty,\theta_{i+1}]$  if
the following conditions are fulfilled:
\begin{enumerate}
\item[(i)]   z(t)  is continuable to $t=\zeta_i;$
\item[(ii)] the derivative $z'(t)$ exists at each point $t \in (- \infty,\theta_{i+1})$ with the possible exception
of the points $\,\theta_j \in (- \infty,\theta_{i+1}),$   where the  one-sided derivatives exist;
\item[(iii)] equation (\ref{1}) is satisfied with $z(t)$ at each point $t \in (- \infty,\theta_{i+1})\backslash \{\theta_i\},$ and at the points $\theta_j \in (- \infty,\theta_{i+1})$
it  holds for the right derivative of $z(t).$ 
\end{enumerate}
\label{defn+3}
\end{definition}
We shall also use  the following definition, which is a version of a definition from \cite{p1}, modified  for our general case.
\begin{definition}   A function $z(t)$ is a solution of (\ref{1}) on $\mathbb R$ if:
\begin{enumerate}
\item[(i)] $z(t)$ is continuous on $\mathbb R;$
\item[(ii)] the derivative $z'(t)$ exists at each point $t \in \mathbb R$ with the possible exception
of the points $\theta_i , i \in \mathbb Z,$  where the one-sided derivatives exist;
\item[(iii)] equation (\ref{1})  is satisfied for $z(t)$ on each interval $(\theta_i, \theta_{i+1}), i \in \mathbb Z,$ and it holds for the  right derivative of $z(t)$ at the points $\theta_i , i \in \mathbb Z.$ 
\end{enumerate}
\label{defn99}
\end{definition}
\begin{definition} 
A set $\Sigma$  in the $(t,z)-$ space is said to be  an  integral set of system (\ref{1}) if any solution
$z(t) = z(t,t_0,z_0), z(t_0) = z_0,$ with $(t_0,z_0) \in \Sigma,$ has the property that $(t,z(t)) \in \Sigma, t \in \mathbb R.$ In other words, for every $(t_0,z_0) \in \Sigma$ the solution $z(t) = z(t,t_0,z_0), z(t_0) = z_0,$ is continuable on  $\mathbb R$ and 
$(t,z(t)) \in \Sigma, t \in \mathbb R.$
\label{defn54}
\end{definition}
\begin{definition} 
A set $\Sigma$  in the $(t,z)-$ space is said to be  a local integral set of system (\ref{1}) if for every  $(t_0,z_0) \in \Sigma$ there exists $\epsilon >0,\, \epsilon = \epsilon(t_0,z_0),$ such 
that if  $z(t) = z(t,t_0,z_0)$  is a solution of  (\ref{1}) and  $|t-t_0| < \epsilon$ then $(t,z(t)) \in \Sigma.$ 
\label{defn55}
\end{definition}

\subsection{The existence and uniqueness of solutions on $\mathbb R$}
In what follows we use the uniform norm $||T|| = \sup\{||Tz|| | ||z||\leq 1\}$ for matrices.

It is known that there exists a constant $\Omega>0$ such that 
$||{\rm e}^{A(t-s)}|| \leq {\rm e}^{\Omega|t-s|}, \\t,s \in \mathbb R.$
Hence, one can show that
\begin{eqnarray*}
&& ||e^{A(t-s)}|| \geq  {\rm e}^{-\Omega|t-s|}, t,s \in \mathbb R.
\label{w2}
\end{eqnarray*}

The last two  inequalities imply  the following,  very simple but useful in  what follows,  estimates
\begin{eqnarray*}
&& ||e^{A(t-s)}|| \leq M, \quad ||e^{A(t-s)}|| \geq m,
\label{w3}
\end{eqnarray*}
if $|t-s| \leq \theta,$ where $M =  {\rm e}^{\Omega \theta}, m =  {\rm e}^{-\Omega \theta}.$

From now on  we make the assumption:
\begin{itemize}
	\item [(C5)]$Ml\theta {\rm e}^{Ml\theta}<1,\\2Ml\theta<1,\\  M^2l\theta \Big \{
	\frac{ Ml\theta {\rm e}^{Ml\theta}+1}{1- Ml\theta {\rm e}^{Ml\theta}} + 
	 Ml\theta {\rm e}^{Ml\theta}\Big \} < m.$
\end{itemize}
\begin{theorem} 
Assume that conditions $(C1)-(C3),$ and $(C5)$ are fulfilled.
Then for  every $(t_0,z_0) \in \mathbb R \times \mathbb R^n$ there
exists a solution $z(t)=z(t,t_0,z_0)$ of  (\ref{1}) which is defined on  
$\mathbb R$ and is unique.
\label{uay}
\end{theorem}
{\it Proof.} 
{\it The existence of the solution.} Let us consider only backward continuation, since forward continuation can be investigated in a similar manner. Theorem \ref{gu-gu} implies that it is sufficient to consider the continuation of a solution $z(t) = z(t,\theta_i,z(\theta_i,t_0,z_0))$ from $\theta_i$ to $\zeta_{i-1},$ for all $ i \in \mathbb Z.$ We have that 
\begin{eqnarray*}
z(t) = {\rm e}^{A(t- \theta_i)}z(\theta_i)
 + \int_{\theta_i}^{t}{\rm e}^{A(t-s)} f(s,z(s),z(\zeta_{i-1}))ds
\label{sec}	
\end{eqnarray*}
on $ [\zeta_{i-1}, \theta_i].$ 

Define a norm $||z(t)||_0 = \max_{[\zeta_{i-1}, \theta_i]}||z(t)||,$ and take 
$  z_0(t) = e^{A(t- \theta_i)}z(\theta_i).$ Define a sequence 
\begin{eqnarray*}
z_{m+1}(t) = {\rm e}^{A(t- \theta_i)}z(\theta_i) + \int_{\theta_i}^{t}{\rm e}^{A(t-s)} f(s,z_m(s),z_m(\zeta_{i-1}))ds, \, m \geq 0.
\label{sec1}	
\end{eqnarray*}
The last expression implies that 
\begin{eqnarray*}
&& ||	z_{m+1}(t)- 	z_{m}(t)||_0 \leq [2Ml\theta]^{m+1} M ||z(\theta_i)||.
\label{trap}
\end{eqnarray*}
The existence is proved.

{\it The uniqueness of the solution.} Denote  $z_j(t) = z(t,t_0,z_0^j), z_j(t_0) = z_0^j, j =1,2,$  solutions of 
(\ref{1}), where $\theta_i\leq t_0 \leq \theta_{i+1}.$ 
It is sufficient to check  that for every $t \in [\theta_i, \theta_{i+1}], $
$ z_0^1 \not = z_0^2$ implies $ z_1(t)\not = z_2(t).$
We have that 
$$z_1(t)- z_2(t) = {\rm e}^{A(t- \theta_i)} (z_0^2 - z_0^1)  - \int_{t_0}^{t}{\rm e}^{A(t-s)}[ f(s,z_1(s),z_1(\zeta_i)- f(s,z_2(s),z_2(\zeta_i)]ds.$$
Hence, 
$$||z_1(t)- z_2(t)|| \leq M ||z_0^2 - z_0^1|| + M l \theta||z_1(\zeta_i)- z_2(\zeta_i)||
+  Ml |\int_{t_0}^{t} ||z_1(s)-z_2(s)||ds|.$$
The Gronwall-Bellman Lemma yields that
\begin{eqnarray*}
||z_1(t)- z_2(t)|| \leq M (||z_0^2 - z_0^1|| +  l \theta||z_1(\zeta_i)- z_2(\zeta_i)||){\rm e}^{Ml\theta}.	
\label{op}
\end{eqnarray*}
Particularly, 
$$||z_1(\zeta_i)- z_2(\zeta_i)|| \leq M (||z_0^2 - z_0^1|| +  l \theta||z_1(\zeta_i)- z_2(\zeta_i)||){\rm e}^{Ml\theta} .$$
Then, 
$$||z_1(\zeta_i)- z_2(\zeta_i)|| \leq \frac{M}{1-Ml\theta {\rm e}^{Ml\theta}}||z_0^2 - z_0^1||.$$
Hence,
\begin{eqnarray}
||z_1(t)- z_2(t)|| \leq 
M {\rm e}^{Ml\theta}[1+ \frac{ Ml\theta}{1- Ml\theta {\rm e}^{ Ml\theta}}]||z_0^2 - z_0^1||.	
\label{opa}
\end{eqnarray}
Assume on the contrary that there exists $t \in [\theta_i, \theta_{i+1}]$ such that $z_1(t) = z_2(t).$
Then
\begin{eqnarray} 
&&{\rm e}^{A(t- t_0)} (z_0^1 - z_0^2)  = \nonumber\\
&&\int_{t_0}^{t}{\rm e}^{A(t-s)}[ f(s,z_2(s),z_2(\zeta_i))- f(s,z_1(s),z_1(\zeta_i))]ds.
\label{bw1}
\end{eqnarray}
We have that
\begin{eqnarray} 
&& ||{\rm e}^{A(t- t_0)} (z_0^2 - z_0^1)|| \geq m ||z_0^2 - z_0^1||.
\label{bw2}
\end{eqnarray}
Moreover, (\ref{opa}) implies that
\begin{eqnarray} 
&&||\int_{t_0}^{t}{\rm e}^{A(t-s)}[ f(s,z_1(s),z_1(\zeta_i)- 
f(s,z_2(s),z_2(\zeta_i)]ds|| \leq \nonumber\\
&&  M^2l\theta \Big \{
	\frac{Ml\theta {\rm e}^{Ml\theta}+1}{1- Ml\theta {\rm e}^{Ml\theta}} + 
	 Ml\theta {\rm e}^{ Ml\theta}\Big \}||z_0^2 - z_0^1||. 
\label{bw3}
\end{eqnarray}
Finally, one can see that $(C5)$,  (\ref{bw2}) and (\ref{bw3})  contradict   (\ref{bw1}).
The theorem is proved.
\begin{remark}
Inequality (\ref{opa}) implies continuous dependence of solutions of 
(\ref{1}) on the initial value.
\label{rem1}
\end{remark}

\section{The existence of integral surfaces}
\label{sect2}
 Fix a number $\sigma \in \mathbb R$ such that
$\mu <- \sigma <0.$ Clearly, there exist
constants $ K\ge 1$ and $m \in \mathbb N, m < n-k,$  such that
\begin{eqnarray*}
&&||{\rm e}^{B_+t}|| \leq K {\rm e}^{-\sigma t},\,\mbox{and} \quad ||{\rm e}^{-B_- t}|| \leq K (1+ t^{m}),
\label{3}
\end{eqnarray*}
for all $t\in R_+ = [0, \infty).$

Using condition  $(C4)$ one can write  equation
(\ref{1})  as  the following system
\begin{eqnarray}\nonumber
\frac{du}{dt} &=&B_+ u + f_+(t,z(t),z(\beta(t))), \nonumber\\
\frac{dv}{dt} &=&B_-v + f_-(t,z(t),z(\beta(t))), 
\label{e3}
\end{eqnarray}
where $z=(u,v), u \in {R}^k , v \in {R} ^{n-k}, 
(f_+,f_-) =f(t,z(t),z(\beta(t))).$

Fix a number $\alpha, 0< \alpha < \sigma,$ and denote $$ \gamma = \int^{\infty}_0 (1+ t^{m}){\rm e}^{-\alpha t} dt.$$
We shall establish the validity of the  following lemma.
\begin{lemma} Fix $ N \in \mathbb R, N >0,$ and assume that conditions $(C1)-(C3)$
 are valid.  
A continuous function $z(t) = (u,v), ||z(t)|| \leq N {\rm e}^{-\alpha (t-t_0)},
t \geq t_0,\,$ is a solution of (\ref{1}) on $\mathbb R$ if and only if
$z(t)$ is a solution 
on $\mathbb R$ of the following
system of integral equations
\begin{eqnarray}
&& u(t) = {\rm e}^{B_+(t-t_0)}u(t_0) + \int^t_{t_0}e^{B_+(t-s)}f_+(s,z(s),z(\beta(s)))ds,\nonumber\\
&& v(t) = - \int^{\infty}_t {\rm e}^{B_-(t-s)}f_- (s,z(s),z(\beta(s)))ds.
\label{6}
\end{eqnarray}
\label{l1}
\end{lemma}
{\it Proof.}  \rm {\it Necessity.} Assume that $z(t) = (u,v), ||z(t)|| \leq N {\rm e}^{ -\alpha (t-t_0)}, t \in [t_0, \infty),$ is 
a solution of (\ref{1}). Denote

\begin{eqnarray}
&& \phi(t)  = {\rm e}^{B_+(t-t_0)}u(t_0) + \int^t_{t_0}{\rm e}^{B_+(t-s)}f_+(s,z(s),z(\beta(s)))ds,\nonumber\\
&& \psi(t) =  - \int^{\infty}_t{\rm e}^{B_-(t-s)}f_- (s,z(s),z(\beta(s)))ds.
\label{7}
\end {eqnarray}
By straightforward evaluation we can see that the integrals converge, are bounded on $[t_0, \infty),$ and, moreover, 
\begin{eqnarray}
&&||\phi(t)|| \leq K {\rm e} ^{-\sigma(t-t_0)}||u(t_0)|| + \nonumber\\
&& Kl \Big [\frac{N(1+{\rm e}^{\alpha\theta})}{\sigma -\alpha} + 
2\frac{1}{\sigma} \max_{[\theta_i, \theta_{i+1}]} ||z(s)||{\rm e}^{\sigma(|t_0| + \theta)} \
\Big]  e ^{-\alpha(t-t_0)},\nonumber\\
&& ||\psi(t)|| \leq Kl \gamma N(1+{\rm e}^{\alpha\theta}) {\rm e} ^{-\alpha(t-t_0)}.
\label{123} 	
\end{eqnarray}
If $t \not = \theta_i, i \in \mathbb Z,$ then 
\begin{eqnarray*}
&& \phi'(t)  = B_+\phi(t) +  f_+(t,z(t),z(\beta(t))), \nonumber\\
&& \psi'(t) = B_-\psi(t) +  f_-(t,z(t),z(\beta(t))),
\end {eqnarray*}
and 
\begin{eqnarray*}
&& u'(t)  = B_+u(t) +  f_+(t,z(t),z(\beta(t))), \nonumber\\
&& v'(t) = B_-v(t) +  f_-(t,z(t),z(\beta(t))).
\end {eqnarray*}
Hence, 
\begin{eqnarray*}
&& [ \phi(t) -u(t)]'  = B_+[ \phi(t) -u(t)], \nonumber\\
&& [ \psi(t) -v(t)]'  = B_-[ \psi(t) -v(t)].
\end {eqnarray*}
Calculating the limit values at  $\theta_j \in \mathbb Z$ we can find that
$$\phi'(\theta_j \pm 0)  = B_+\phi(\theta_j \pm 0) + f_+(\theta_j \pm 0,z(\theta_j \pm 0),z(\beta(\theta_j \pm 0))),$$
$$u'(\theta_j \pm 0)  = B_+u(\theta_j \pm 0) +  f_+(\theta_j \pm 0,z(\theta_j \pm 0),z(\beta(\theta_j \pm 0))),$$
$$ \psi'(\theta_j \pm 0)  = B_+ \psi(\theta_j \pm 0) + f_(\theta_j \pm 0, z(\theta_j \pm 0), z(\beta(\theta_j \pm 0))),$$
$$v'(\theta_j \pm 0)  = B_+v(\theta_j \pm 0) +  f_-(\theta_j \pm 0,z(\theta_j \pm 0),z(\beta(\theta_j \pm 0))).$$
Consequently,
$$ [ \phi(t) -u(t)]'|_{ t=\theta_j+0} =    [ \phi(t) -u(t)]'|_{ t=\theta_j-0},  [ \psi(t) -v(t)]'|_{ t=\theta_j+0} =    [ \psi(t) -v(t)]'|_{ t=\theta_j-0}.$$
Thus,  $(\phi(t) -u(t),   \psi(t) -v(t)) $ is a continuously differentiable on $\mathbb R$ function  satisfying
$u'(t) = B_+ u(t),  v'(t) = B_-v(t)$ with the initial condition $\phi(t_0) -u(t_0) = 0.$ Assume that  $\psi(t_0) -v(t_0) \not =0.$
 Then $\psi(t) -v(t)$ is not  a decay solution, which contradicts (\ref{123}). Hence, $\phi(t) -u(t)=0,   \psi(t) -v(t)=0$ 
on $\mathbb R.$\\
{\it Sufficiency.}  Suppose that $z(t)$ is a solution of (\ref{6}).  Differentiating 
$z(t)$ in $ t \in (\theta_i, \theta_{i+1}),\, i \in \mathbb Z,$  one can see that the function satisfies (\ref{1}).
Moreover, letting  $ t \rightarrow \theta_i+,$  and remembering  that $z(\beta(t))$
is a right-continuous function, we find that $z(t)$ satisfies (\ref{1}) on $[\theta_i, \theta_{i+1}).$ The Lemma is proved.

Denote $$p= K(1+ {\rm e}^{\alpha \theta})[\frac{1}{\sigma - \alpha} + 
  \gamma].$$
In what follows we mainly use the technique of \cite{pliss1}. See  also \cite{carr,pal1}.
\begin{theorem} Suppose conditions $(C1)-(C5)$ are fulfilled and, moreover,
\begin{eqnarray}
&& 2pl < 1.
\label{10}
\end{eqnarray}
 Then  for arbitrary $\alpha \in (0, \sigma)$  there 
exists a  function
$F(\zeta_i,u),\, i \in \mathbb Z,$  satisfying
\begin{eqnarray}
&& F(\zeta_i,0)=0,\label{tip-tip}\\
&& ||F(\zeta_i,u_1) - F(\zeta_i,u_2)||\leq pKl||u_1-u_2||,
\label{top-top}
\end{eqnarray}

for all $i, u_1, u_2,$  such that   a solution $z(t)$ of  (\ref{1}) with  $z(\zeta_i) = (c, F(\zeta_i,c), c \in \mathbb R^k,$ is defined on $\mathbb R$ and satisfies 
\begin{eqnarray}
&& ||z(t)|| \leq 2K ||c|| {\rm e}^{-\alpha(t-\zeta_i)}, t \geq \zeta_i.
\label{7in}
\end{eqnarray}
\label{+t1}
\end{theorem}
{\it Proof.} Let us consider  system (\ref{6}) and apply the method of successive approximations.
Denote $z_0(t)=(0, 0)^T, z_m = (u_m,v_m)^T, m \in \mathbb N,$ where for $ m \geq 0$
\begin{eqnarray*}
&& u_{m+1}(t) = {\rm e}^{B_+(t-\zeta_i)}c + \int^t_{\zeta_i}{\rm e}^{B_+(t-s)}f_+(s,z_m(s),z_m(\beta(s)))ds,\nonumber\\
&& v_{m+1}(t) = - \int^{\infty}_t {\rm e}^{B_-(t-s)}f_- (s,z_m(s),z_m(\beta(s)))ds.
\label{8}
\end{eqnarray*}
Let us show that 
\begin{eqnarray}
&& ||z_m(t)|| \leq 2K ||c|| {\rm e}^{-\alpha(t-\zeta_i)}, t \leq \zeta_i.
\label{7****}
\end{eqnarray}
Indeed, $z_0$ satisfies the relation. Assume that $z_{m-1}$ satisfies (\ref{7****}).
Then 
\begin{eqnarray}
&& ||u_m(t)|| \leq K {\rm e}^{-\sigma(t-\zeta_i)} ||c|| +
\frac{2K^2l(1+ {\rm e}^{\alpha \theta})}{\sigma - \alpha} {\rm e}^{-\alpha(t-\zeta_i)}||c||, \nonumber\\
&& ||v_m(t)|| \leq  2\gamma K^2 l (1+ {\rm e}^{\alpha \theta})  {\rm e}^{-\alpha(t-\zeta_i)}||c||,
\label{1234}
\end{eqnarray}
and  (\ref{7****}) is valid  provided (\ref{10}) is correct.
Similarly, one can establish  the following inequality 
\begin{eqnarray}
 ||z_{m+1}(t) - z_m(t) || \leq K||c||(2pl)^m {\rm e}^{-\alpha(t-\zeta_i)}.
\label{11}
\end{eqnarray}
The last inequality  and assumption (\ref{10})
 imply that the sequence $z_m$ converges uniformly
for all $c$ and $t\geq \zeta_i.$ Let  $z(t,\zeta_i,c)=(u(t,\zeta_i,c), v(t,\zeta_i,c))$ be the limit function. It is obvious that the function is a solution of (\ref{6}). By Lemma \ref{l1} $z(t,\zeta_i,c)$ is a solution of  (\ref{1}), too.
  Taking $t_0=\zeta_i$ in (\ref{6}) we  have that
$$u(\zeta_i,\zeta_i,c) = c,$$$$ v(\zeta_i,\zeta_i,c) = - \int^{\infty}_{\zeta_i} {\rm e}^{B_-(t-s)}f_- (s,z(s,\zeta_i,c),z(\beta(s),\zeta_i,c)))ds.$$
Denote $F(\zeta_i,c) = v(\zeta_i,\zeta_i,c).$ Since 
\begin{eqnarray}
&& ||v_m(t,\zeta_i,c_1) - v_m(t,\zeta_i,c_2) || \leq pKl ||c_1-c_2||, \,m \geq 1,
\label{+11}
\end{eqnarray}
inequality (\ref{top-top}) is valid.
The theorem is proved.

For every $i \in \mathbb Z$ consider a set $\Psi_i$ of continuous on $\mathbb R$ functions   such that
if $\psi \in \Psi_i$ then there exists a positive constant $K_{\psi},$ satisfying 
$||\psi(t)|| \le K_{\psi}{\rm e}^{-\alpha(t-\zeta_i)}, \zeta_i \le t,$ where constant $\alpha$ is defined for Theorem
\ref{+t1}.
\begin{lemma} For every  $\zeta_i, i \in \mathbb Z, c \in \mathbb R^k,$ the system 
\begin{eqnarray*}
&& u(t) = {\rm e}^{B_+(t-t_0)}c + \int^t_{\zeta_i}e^{B_+(t-s)}f_+(s,z(s),z(\beta(s)))ds,\nonumber\\
&& v(t) = - \int^{\infty}_t {\rm e}^{B_-(t-s)}f_- (s,z(s),z(\beta(s)))ds.
\label{6##}
\end{eqnarray*}
has only one solution from  $\Psi_i.$
\label{lemma2}
\end{lemma}
{\it Proof.} \rm If  $z_1$ and $z_2$ are two solutions of the system bounded on $[\zeta_i,\infty),$   then by straightforward evaluation one can  show that
$$ \sup_{[\zeta_i,\infty)} ||  z_1 -  z_2|| \leq 2pl \sup_{[\zeta_i,\infty)} ||z_1 -  z_2||.$$
Hence, in view of (\ref{10})  the lemma is proved.

Let us denote by $S_i^+$ the set of all points from  the $(t,z)-$ space  $(z=(u,v))$ such that
$t=\zeta_i, \, v=F(\zeta_i,u).$ 
\begin{lemma} If $(\zeta_i, z_0) \not \in S_i^+,$ then the solution $z(t,\zeta_i,z_0)$ of (\ref{1})   is  not from  $\Psi_i.$
\label{l3}
\end{lemma} 
{\it Proof.} \rm Assume on the contrary  that  $z(t)\in \Psi_i, z(t) = z(t,\zeta_i,z_0) = (u,v),$ is a solution of (\ref{1}) and $(\zeta_i,z_0) \not \in S_i^+.$ It is obvious that 
\begin{eqnarray}
&& u(t)  = U(t,\zeta_i)u(\zeta_i) + \int^t_{\zeta_i}U(t,s)f_+(s,z(s),z(\beta(s)))ds,\nonumber\\
&& v(t) = V(t,\zeta_i) \hat{\kappa}  - \int^{\infty}_t V(t,s)f_- (s,z(s),z(\beta(s)))ds,
\label{13}
\end {eqnarray}
where
$$\hat{\kappa} = v(\zeta_i) + \int^{\infty}_{\zeta_i} V(t,s)g_- (s,z(s),z(\beta(s)))ds,$$
and the improper integral converges and is bounded on $[\zeta_i,\infty).$
But condition  $(C3)$ on eigenvalues of matrix $B_-$ imply  that 
$||V(t,\zeta_i) \hat{\kappa}|| \to 0$ as $t \to \infty$ if only
$\hat{\kappa} =0.$ By Lemma \ref{l1}  $z(t)$ satisfies (\ref{6}) with $t_0 = \zeta_i.$ The contradiction proves the lemma.

Let $S^+$ be the set of all points from the $(t,z)-$ space $(z=(u,v))$ such that either 
$(t,z)\in S_i^+$ for some $i \in \mathbb Z,$ or there exist \, $\zeta_i,     \zeta_i < t, c \in \mathbb R^k,$ such that $(\zeta_i,c) \in S_i^+ \, {\rm and}\, z = z(t, \zeta_i,c).$
\begin{theorem}  $S^+$ is an  invariant set.
\label{t3}
\end{theorem}
{\it Proof.}  Assume that $(\zeta_i,u_0,v_0) \in S_i^+, z_0 = (u_0,v_0).$ We show that if
$z(t) = z(t, \zeta_i,z_0),$ then $(\theta_j,z(\theta_j)) \in S_j^+$ for  all $j \geq i.$
Indeed, if $t \geq \theta_j$ then  $||z(t)|| \leq  (K+ \epsilon) ||u_0||  {\rm e}^{-\alpha(t-\theta_j)}{\rm e}^{-\alpha(\theta_j-\zeta_i)}.$
Lemma \ref{l1} implies that the  point $(\theta_j,z(\theta_j))$ satisfies the equation $v=F(t,u).$
If $(\theta,\xi)\in S^+\backslash \cup_{i \in \mathbb Z}S_i^+,$ then by the definition and the previous  part of the proof $z(t,\theta,\xi) \in S^+$ for all $t \geq \theta.$ 
Assume that $(\zeta_i,z_0) \in S_i^+,$ and  denote $z(t) = z(t, \zeta_i,z_0).$ Lemma \ref{l3} implies that $(z(\theta_{i-1}), \theta_{i-1}) \in S_j^+.$ 
The theorem is proved.

On the basis of  Theorem \ref{t3}, Lemmas \ref{l1} and \ref{l3} we can conclude that there exists an invariant surface 
$S^+$ of equation (\ref{1}), such that every solution  starting at $S^+$ tends to  zero as  $t\rightarrow \infty.$

Denote $z(t,r,c) = (u(t,r,c), v(t,r,c)), t,r \in \mathbb R, c \in \mathbb R^k,$ a solution
of (\ref{1}) such that $u(r,r,c)=c.$ From the discussion above it can be seen  that
surface $S^+$ contains solutions which satisfy the equation $v=F(t,u), (t,u) \in \mathbb R \times \mathbb R^k,$ where
\begin{eqnarray}
F(t,u) = - \int^{\infty}_t {\rm e}^{B_-(t-s)}f_- (s,z(s,t,u),z(\beta(s),t,u))ds.	
\label{man}
\end{eqnarray}
It is obvious that $F(t,u)$ is a  function continuous in both arguments.

\begin{theorem} Suppose conditions $(C1)-(C5)$ are fulfilled. 
Then  for an arbitrarily   small positive  $\tilde \alpha$ and a sufficiently small Lipschitz
constant $l$ there 
exists  a function
$G(\zeta_i,v), i \in \mathbb Z,$ from $\mathbb R^{n-k}$ to $\mathbb R^{k},$ satisfying
\begin{eqnarray}
&& G(\zeta_i,0)=0,\\
&& ||G(\zeta_i,d_1) - G(\zeta_i,d_2)|| \leq   P l ||d_1-d_2||
\label{in^}
\end{eqnarray}
for all $d_1, d_2,$ such that  a solution $z(t)$  of  (\ref{1}) with 
$z(\zeta_i) = ( G(\zeta_i,v_0),v_0), v_0 \in \mathbb R^{n-k},$ is defined on $\mathbb R$ and satisfies
\begin{eqnarray}
&& ||z(t)|| \leq D ||v_0||{\rm e}^{-\tilde \alpha(t-\zeta_i)},\, t \leq \zeta_i,
\label{7"**}
\end{eqnarray}
where $P, D >0$ are  constant.
\label{t2}
\end{theorem}
{\it Proof.}  Let us denote $\kappa =\frac{\sigma}{2}$ and $\eta(t) = z(t) {\rm e}^{\kappa t}.$ Then system (\ref{1}) is transformed into the equation
\begin{eqnarray}
\frac{d\xi}{dt} &=&(B_++\kappa I)\xi + g_+(t,\eta(t),\eta(\beta(t))), \nonumber\\
\frac{d\omega}{dt} &=&(B_-+\kappa I)\zeta + g_-(t,\eta(t),\eta(\beta(t))), 
\label{=6}
\end{eqnarray}
where $\eta = (\xi, \omega), I$ is an identity matrix, $\eta(\beta(t)) =  z(\beta(t)) {\rm e}^{-\kappa\beta(t)},$ and $
g(t,z,y) = (g_+,g_-) =  {\rm e}^{\kappa t} f(t, z {\rm e}^{-\kappa t}, y {\rm e}^{-\kappa \beta(t}).$
It is easy to see that the function $g(t,z,y)$ satisfies the Lipschitz condition in $z, y$ with a constant $l {\rm e}^{\kappa \theta},$ and the 
eigenvalues of the matrices $B_++\kappa I$ and $B_-+\kappa I$ have negative and positive real parts, respectively, such that
$\mu + \kappa = \max_{j=\overline{1,k}} \Re \lambda_j(B_++\kappa I)< -\sigma + \kappa <0,$ 
and $\min_{j=\overline{k+1,n}} \Re \lambda_j (B_-+\kappa I)=\kappa>0.$
Fix a positive number $\bar \kappa < \min\{\sigma-\kappa, \kappa\} = \kappa.$
There exists a positive number $\bar K$ such that 
\begin{eqnarray*}
&& ||{\rm e}^{(B_++\kappa I)(t-s)}|| \leq \bar K {\rm e}^{-\bar \kappa(t-s)}, t\geq s\nonumber\\
&& ||{\rm e}^{(B_- +\kappa I)(t-s)}|| \leq \bar K {\rm e}^{\bar \kappa(t-s)}, t\leq s.
\label{eval}
\end{eqnarray*}
To continue the proof we  need the following two assertions which can be proved 
similarly to Lemma \ref{l1} and  Theorem \ref{+t1}.
\begin{lemma} Fix $ N \in \mathbb R, N >0,$ and assume that conditions $(C1)-(C3)$
are valid.  
A continuous function $\eta(t) = (\xi,\omega), ||\eta(t)|| \leq N {\rm e}^{\tilde \alpha (t-t_0)}, 0<\tilde \alpha<\bar \kappa, t \leq t_0,\theta_j < t_0 \leq \theta_{j+1},$ is a solution of (\ref{=6}) on $(-\infty,t_0]$ if and only if  $\eta(t)$ is a solution of the following
system of integral equations
\begin{eqnarray}
&& \xi(t) =  \int_{-\infty}^t {\rm e}^{(B_+ + \kappa I)(t-s)}g_+(s,\eta(s),\eta(\beta(s)))ds,\nonumber\\
&& \omega(t) = {\rm e}^{(B_-+ \kappa I)(t-t_0)}\omega(t_0) + \int^{t}_{t_0} {\rm e}^{(B_-+ \kappa I)(t-s)}g_- (s,\eta(s),\eta(\beta(s)))ds.
\label{??6}
\end{eqnarray}
\label{+l3}
\end{lemma}
\begin{lemma} Suppose conditions $(C1)-(C5)$  are fulfilled. Then  for an arbitrary $\alpha_1 \in (0, \bar \kappa)$
and a sufficiently small Lipschitz
constant $l$ there 
exists  a function
$\bar G(\zeta_i,u), i \in \mathbb Z,$ satisfying
\begin{eqnarray}
&& \bar G(\zeta_i,0)=0,\\
&& ||\bar G(\zeta_i,d_1) - \bar G(\zeta_i,d_2)||  \leq Pl||d_1-d_2||,
\label{()in}
\end{eqnarray}
where $P$ is a positive constant, and  such that $\xi_0=\bar G(\zeta_i,\omega_0)$ defines  a solution $\eta(t)$ of  (\ref{1}) with $\eta(\zeta_i) = (\bar G(\zeta_i,\omega_0),\omega_0)$ and 
\begin{eqnarray}
&& ||\eta(t)|| \leq  2\bar K ||\omega_0|| {\rm e}^{\alpha_1(t-\zeta_i)}, t \leq \zeta_i.
\label{15}
\end{eqnarray}
\label{l4}
\end{lemma}
Similarly to (\ref{man}) one can show that if $\eta(t,r,c) = (\xi,\omega)$ is a solution of
(\ref{1}) such that $\omega(r) = c,$ then 
\begin{eqnarray}
&& \bar G(t, \omega) =   \int_{-\infty}^t {\rm e}^{(B_+ + \kappa I)(t-s)}g_+(s,\eta(s,t,\omega),\eta(\beta(s),t,\omega))ds.
\label{???6}
\end{eqnarray}

Let us now finish  the proof of Theorem \ref{t2}. Applying the inverse transformation $z(t) = \eta(t) {\rm e}^{- \kappa t}$ we
can define a new function $G(\zeta_i,v) ={\rm e}^{- \kappa t} \bar G (\zeta_i, v {\rm e}^{\kappa t})$ and  check that
\begin{eqnarray}
&& ||G(\zeta_i,v_1) -G(\zeta_i,v_2) || \leq  Pl ||v_1-v_2||,
\label{tak}
\end{eqnarray}
 and 
$$||z(t)|| \leq  2\bar K ||v_0|| {\rm e}^{(\alpha_1-\kappa)(t-\zeta_i)}, t \leq \zeta_i,$$
if $u_0 = G(\zeta_i, v_0).$  If we denote now $D = 2 \bar K$ and  choose $\bar \kappa$ sufficiently close to $\kappa$ then we can take $\alpha_1 = \kappa - \tilde \alpha > 0$ such that the last inequality implies
(\ref{7"**}). The theorem is proved.

Using the equation $\xi = \bar G(t, \omega)$  we can define, similarly to $S^+$, an integral surface $\bar S_0$ such that every solution of (\ref{=6}) starting on  $\bar S_0$ tends to the origin as $t \rightarrow -\infty.$ Then an integral  set $S_0$ for (\ref{1}) can be  defined by the equation $u = G(t, v).$

\section{The stability of the zero solution}

We shall need the following definitions.
\begin{definition} The trivial solution  of (\ref{1}) is stable, if for any $\epsilon >0$ and any $t_0 \in \mathbb R,$ there exists a $\delta(t_0,\epsilon) > 0$  such that if 
$||z_0|| < \delta(t_0, \epsilon),$  then $||z(t,t_0,z_0)|| < \epsilon$ for all $t \ge t_0.$
If the $\delta$ above is independent of $t_0$ then the zero solution is uniformly stable.
\label{stable}
\end{definition}
\begin{definition} The zero solution of  (\ref{1}) is asymptotically stable, if it is stable and if there exists a $\delta_0(t_0) > 0$  such that if 
$||z_0|| < \delta_0(t_0),$  then $z(t,t_0,z_0) \to 0,$ as $t \to \infty.$
\label{asymstable}
\end{definition}
\begin{definition} The zero solution of  (\ref{1}) is uniformly asymptotically stable, if it is 
uniformly stable and there is a $\kappa_0 > 0$  such that  for any $t_0 \in \mathbb R,$ there
exists  a $T(\epsilon) >0,$ independent of $t_0,$ such that $||z(t,t_0,z_0)|| < \epsilon$ for all $t \ge t_0 + T(\epsilon)$ whenever $||z_0|| < \kappa_0.$ 
\label{uasymstable}
\end{definition}
\begin{definition} The zero solution of  (\ref{1}) is exponentially stable if there exists an $\alpha > 0,$ and  for every $\epsilon >0$ and $t_0$  there exists a $\delta(\epsilon,t_0) >0,$ such that 
$$||z(t,t_0,z_0)|| \le  \epsilon {\rm e}^{-\alpha(t-t_0)}$$
for all $t \ge t_0,$ whenever $||z_0|| < \delta.$  If the $\delta$ above is independent of $t_0$ then the zero solution is uniformly exponentially stable.
 
\label{expstable}
\end{definition}

System (\ref{1}) is  an equation with a deviating argument, but one can easily see that  Definitions \ref{stable} - \ref{expstable} coincide with  the definitions of stability in the Lyapunov sense for ordinary differential equations \cite{hart, mil}. They do not involve the concept of initial interval for an initial value problem. This phenomenon must not surprise us,  as  the  right side of (\ref{1}) depends only on one "delayed" value of a solution at $t=\zeta_i$ if $\theta_i \le t < \theta_{i+1}, i \in \mathbb Z.$ For EPCA where argument is delayed  \cite{cwm,w} the stability is investigated   with $t_0 = 0.$ 
Continuous dependence on the initial value provided by (\ref{opa}) helps us to investigate stability assuming that the initial moment $t_0$ can be an arbitrary real number.

Theorem \ref{+t1} considered with $k=n$ and inequality (\ref{opa}) imply that the 
following assertion is valid.

\begin{theorem} Suppose that conditions $(C1),(C2)$ and $(C5)$ are fulfilled, and all eigenvalues
of matrix $A$ have negative real parts.
Then the zero solution of (\ref{1}) is uniformly exponentially stable if the Lipschitz constant $l$ is sufficiently small.
\end{theorem}

Comparing Definitions \ref{uasymstable} and \ref{expstable} we can conclude that if  the zero solution
is  uniformly exponentially stable then it is uniformly asymptotically stable.


\section{The stability of the integral surface $S_0$}

\begin{theorem} 
If the Lipschitz constant $l$ is sufficiently small then for every  solution $z(t) = (u,v)$ of (\ref{1}) there exists  a solution $\mu(t) = (\phi, \psi)$ on   
 $S_0$ such that
\begin{eqnarray}
&& ||u(t) - \phi(t)|| \leq 2 K||u(\zeta_i) - \phi(\zeta_i)|| {\rm e}^{ -\alpha(t-\zeta_i)},\nonumber\\
&&  ||v(t) - \psi(t)|| \leq K ||v(\zeta_i) - \psi(\zeta_i)|| {\rm e}^{-\alpha(t-\zeta_i)}, 
\zeta_i \le t,
\label{23}
\end {eqnarray}
where  $\alpha$ is the coefficient defined for Theorem \ref{+t1}.
\label{t65}
\end{theorem}
{\it Proof.} Fix a solution $z(t,\zeta_i,z_0)$ of (\ref{1}). Denote by 
$\mu(t, \zeta_i,d) = (\phi, \psi)$ a solution of  (\ref{1}) such that $\psi(\zeta_i,\zeta_i,d) = d, \phi(\zeta_i,\zeta_i,d) = G(\zeta_i,d).$
Let us carry out  the transformation
\begin{eqnarray}
&& X(t) = u - \phi(t), Y(t) = v - \psi(t), 
\label{tr}
\end{eqnarray}
in system  (\ref{1}) and denote $Z = (X,Y).$  The transformed equation has the form
\begin{eqnarray}
&& \frac{dX}{dt} = B_+ X + Q_+(t,Z(t), Z(\beta(t)),d), \nonumber\\
&& \frac{dY}{dt} = B_-Y + Q_-(t,Z(t), Z(\beta(t)),d),
\label{treqn}
\end{eqnarray}
where $Q(t,X,Y,d) = (Q_+,Q_-) = f(t,z(t),z(\beta(t))) - f(t,\mu(t),\mu(\beta(t))).$ One can see that 
$Q$ satisfies the Lipschitz condition with the same  constant $l.$ By Theorem \ref{+t1}  there exists a function
$\tilde{F}$ such that the  equation $Y = \tilde{F}(\zeta_i,X,d)$ defines a set for (\ref{treqn}) which satisfies, according to (\ref{tip-tip}) and (\ref{top-top}), the following properties
\begin{eqnarray}
&& \tilde{F}(\zeta_i,0,d) = 0, \nonumber\\
&& || \tilde{F}(\zeta_i,X_1,d) -  \tilde{F}(\zeta_i,X_2,d)|| \leq pKl ||X_1-X_2||.
\label{F}
\end{eqnarray}
Using (\ref{10}) and  formulas similar to (\ref{1234})   one can see that every solution \\
$Z(t), Z(\zeta_i) = (X_0, \tilde{F}(\zeta_i, X_0)),$  satisfies 
\begin{eqnarray}
&& ||X(t)|| \leq K {\rm e}^{-\sigma(t-\zeta_i)}||X_0|| +
\frac{2K^2l(1+ {\rm e}^{\alpha \theta})}{\sigma - \alpha} {\rm e}^{-\alpha(t-\zeta_i)}||X_0||, \nonumber\\
&& ||Y(t)|| \leq   2\gamma K^2 l (1+ {\rm e}^{\alpha \theta})  {\rm e}^{-\alpha(t-\zeta_i)}||X_0||.
\label{Imp}
\end{eqnarray}
Let us show that there exist $X_0$ and $d$ such that for solutions $z(t)$ and $(X(t), Y(t))$ of systems
(\ref{1}) and (\ref{treqn}), respectively, 
\begin{eqnarray*}
&& X(t) = u(t,\zeta_i,z_0) - \phi(t), Y(t) = v(t,\zeta_i,z_0) - \psi(t).
\label{tr**}
\end{eqnarray*}
The last equalities for $t = \zeta_i$ have the form 

\begin{eqnarray}
&& X_0 = u_0 - G(\zeta_i, d), \tilde{F}(\zeta_i, X_0,d)  = v_0 - d.
\label{tr***}
\end{eqnarray}
Let us consider the system as an equation with respect to  $X_0$ and $ d.$      We shall show that  it has a solution 
for every pair $(u_0, v_0).$
Equation (\ref{tr***}) implies that

\begin{eqnarray}
&& d = v_0 -  \tilde{F}(\zeta_i, u_0-  G(\zeta_i, d),d).
\label{trac}
\end{eqnarray}
Applying properties (\ref{F}) of the function $\tilde{F}$ and equality (\ref{trac}) we can write that
$$||d-v_0|| \leq pKl||u_0 - G(\zeta_i,d)||.$$
Since the function $G$ satisfies the Lipschitz condition, using
$$||d-v_0|| \leq pKl||u_0 - G(\zeta_i,v_0)|| +  pKl||G(\zeta_i,d) - G(\zeta_i,v_0)||,$$ one can show  that

\begin{eqnarray}
&& ||d-v_0|| \leq \frac{pKl}{1- pPKl^2} ||u_0 - G(\zeta_i,v_0)||.
\label{harsh}
\end{eqnarray}

We assume that $1- pPKl^2 >0, pKl(1+Pl)  \le 1,$ and will consider  the ball 
$\hat B =\{d: ||d-v_0|| \leq ||u_0 -  G(\zeta_i,v_0)||\}.$ Inequality (\ref{harsh}) implies that
 (\ref{trac}) transforms $\hat B$ into itself, and by Brauer's theorem there exists a fixed point of the transformation. Denote the point by $\bar d.$ Substituting $\bar d$ into the first equation of
(\ref{tr***}) one can obtain the value $\bar X_0.$ The pair $(\bar X_0, \bar d)$ satisfies system 
(\ref{tr***}). Now, applying (\ref{Imp}), (\ref{10}) and the theorem of existence and uniqueness we can complete the proof of the theorem.

We shall introduce a notion of stability for an integral set \cite{carr,pliss}.
Denote by $ M \subset \mathbb R \times \mathbb R^n$ an integral surface of (\ref{1}) and by 
$d(z, M)$ the distance between a point $z \in \mathbb R^n$ and the set $M.$
\begin{definition} $M$ is a stable integral surface of  (\ref{1}), if for any  $\epsilon >0,$ there exists a number $\delta>0, \delta = \delta(\epsilon, t_0),$  such that if $d(z_0, M(t_0)) < \delta,$
then \\$d(z(t,t_0,z_0), M(t)) < \epsilon$ for all $t \geq t_0.$ 
\end{definition}
\begin{definition} A stable integral surface $M$ is stable in large, if every solution 
of (\ref{1})  approaches $M$ as $t \rightarrow \infty.$ 
\end{definition}

Theorem \ref{t65}  implies that the surface $S_0$ is stable, and, moreover, it is stable in large.


\section{The reduction principle}
The following conditions are needed in this part of the paper.
\begin{itemize}
\item[(C6)] The function $f(t,z,w)$ is uniformly continuously differentiable in $z,w$ for all
$t, z, w,$  and $$\frac{\partial f(t,0,0)}{\partial z} =0, \frac{\partial f(t,0,0)}{\partial w} =0.$$ 
\item[(C7)]  If we denote by $\lambda_j, j= \overline{1,n},$ the eigenvalues of  matrix $A,$ then there exists a positive integer $k$ such that $\mu = \max_{j=\overline{1,k}} \Re \lambda_j <0,$ 
and $\Re \lambda_j =0,j=\overline{k+1,n},$
where $\Re \lambda_j$ denotes the real part of the eigenvalue $\lambda_j$ of  matrix $A.$

\end{itemize}

Denote 
$$T(h) = \{(t,z) \in \mathbb R \times \mathbb R^n: ||z|| < h\}$$
for a fixed number $h >0.$ Assume that  $\epsilon_0>0$ is sufficiently  small for  the Lipschitz constant $l,$
provided by $(C6),$ to satisfy all conditions of Theorem \ref{t2} in $T(\epsilon_0).$

Denote $\epsilon_1 = \frac{\epsilon_0}{2\bar K},$ where $\bar K$ is the constant from 
(\ref{15}).

By Lemma  \ref{l4}  there exists a local integral manifold 
of (\ref{=6}) in  $T(\epsilon_1)$  such that a solution starting on the manifold  is continuable 
to $-\infty,$ and is exponentially decaying.

Using the inverse transformation $z= \eta{\rm e}^{-\kappa t}$ one can obtain a local integral manifold 
of (\ref{1}) in $T(\epsilon_1)$  given by equation $u= G(t,v).$ Solutions of (\ref{1}) on the manifold
are not necessarily continuable to $-\infty$  in $T(\epsilon_1).$  
For the function $G$ condition (\ref{tak}) is true  and $G(t,0) = 0, t \in \mathbb R.$
On the local manifold solutions of  (\ref{1})
 satisfy the following system
\begin{eqnarray}
\frac{dv}{dt} &=&B_-v + f_-(t,(G(t,v(t)),v(t)),(G(t,v(\beta(t))),v(\beta(t))). 
\label{3*}
\end{eqnarray}
We can see that   the function $f_-(t,(G(t,v),v),(G(t,\bar v), \bar v))  $  satisfies the Lipschitz condition in $v, \bar v$ with the constant 
$l(1+ Pl).$ 

\begin{theorem} Assume that conditions $(C1)-(C2),(C4)-(C7)$ are fulfilled.
The  trivial solution of (\ref{1}) is  stable, asymptotically stable or unstable in Lyapunov sense, if the trivial solution  of (\ref{3*}) is stable, asymptotically stable or unstable, respectively.
\label{9898}
\end{theorem}
{\it Proof.} Consider system  (\ref{1})  in $T(\epsilon_1).$ We assume, additionally, that $\epsilon_0$  is sufficiently small such that conditions of Theorem   \ref{t65}
are valid in  $T(\epsilon_0),$  and, moreover, 
\begin{eqnarray}
&& 1+Pl \le 2.
\label{kryak}
\end{eqnarray} Suppose that the zero solution of  (\ref{3*}) is stable in  the sense of Definition \ref{stable}. Fix an $\epsilon>0.$ Without loss of 
generality we assume that $\epsilon < \epsilon_1.$ 

In  view of Remark \ref{rem1} we can  assume that $t_0=\zeta_i$ for 
some fixed $i \in \mathbb Z.$ 
Fix a positive number $\nu$ such that the 
inequality 
\begin{eqnarray}
&& 2\nu(1+Pl) < 1
\label{lop-top}
\end{eqnarray} is true. 
The stability implies  the existence of $\delta >0, 0< 2\delta < \epsilon,$ such that
if $d \in \mathbb R^{n-k}, ||d|| < 2\delta,$ then the solution $v = \psi(t, \zeta_i,d)$ of (\ref{3*}) satisfies the inequality
\begin{eqnarray}
&& ||\psi(t, \zeta_i,d|| < \nu \epsilon,\, \zeta_i \le t.
\label{rp1}
\end{eqnarray}
Let $u_0$ and $v_0$ be arbitrary vectors satisfying 
$||u_0|| + ||v_0|| < \delta.$ 
Denote $z(t) = z(t,\zeta_i,z_0), z(\zeta_i) = z_0, z_0 = (u_0,v_0),$
a solution of (\ref{1}).
Further we shall follow the proof of Theorem \ref{t65} specifying it for the local case.
Let $\mu(t) =\mu(t, \zeta_i,d) = (\phi, \psi)$ be a solution 
of  (\ref{1}) such that $\psi(\zeta_i, \zeta_i,d) = d, \phi(\zeta_i, \zeta_i,d) = G(\zeta_i,d)$
and $\psi(\zeta_i, \zeta_i,d)$ satisfies (\ref{rp1}).
Applying (\ref{rp1}) and the Lipschitz condition on $G$ we have that
$||\phi(t, \zeta_i,d)|| \le  Pl \nu \epsilon, \zeta_i \le t.$
Then $||\mu((t, \zeta_i,d)|| \le (1+ Pl)\nu \epsilon, \zeta_i \ge t.$
Finally using (\ref{lop-top}) we can write  that
\begin{eqnarray}
&& ||\mu(t)|| < \frac{1}{2}\epsilon.
\label{mutu}
\end{eqnarray}
Applying transformation (\ref{tr}) we obtain equation (\ref{treqn}).
From (\ref{mutu})  it follows that (\ref{tr}) transforms neighborhood $T( \frac{\epsilon}{2})$ for (\ref{treqn}) into neighborhood $T(\epsilon)$ for (\ref{1}).
So, the conditions set  by Theorem \ref{t65} for the  coefficient $l$  are valid  if (\ref{treqn})
is considered  
in $\frac{\epsilon}{2}-$neighborhood  of  $ X= 0, Y = 0, t \in \mathbb R.$

Now, if we assume that
\begin{eqnarray}
&&||X(\zeta_i)|| < \frac{\epsilon}{2K(1+ 2pl)},
\label{chok-chok}
\end{eqnarray}
then similarly to the sequence $(u_m,v_m)$ in Theorem \ref{+t1} we can construct  a sequence 
$Z_m = (X_m,Y_m), m \ge 0,$ such that $(X_0,Y_0) = (0,0)^T,$ 

\begin{eqnarray*}
&& X_{m+1}(t) = {\rm e}^{B_+(t-\zeta_i)}X(\zeta_i) + \int^t_{\zeta_i}{\rm e}^{B_+(t-s)}Q_+(s,Z_m(s),Z_m(\beta(s)))ds,\nonumber\\
&& Y_{m+1}(t) = - \int^{\infty}_t {\rm e}^{B_-(t-s)}Q_- (s,Z_m(s),Z_m(\beta(s)))ds,
\label{81}
\end{eqnarray*}
\begin{eqnarray}
&& ||X_m(t)|| \leq K {\rm e}^{-\sigma(t-\zeta_i)} ||X(\zeta_i)|| +
\frac{2K^2l(1+ {\rm e}^{\alpha \theta})}{\sigma - \alpha} {\rm e}^{-\alpha(t-\zeta_i)}||X(\zeta_i)||, \nonumber\\
&& ||Y_m(t)|| \leq  2\gamma K^2 l (1+ {\rm e}^{\alpha \theta})  {\rm e}^{-\alpha(t-\zeta_i)}||X(\zeta_i)||,
\label{123**}
\end{eqnarray}
and, hence,  
\begin{eqnarray*}
&& ||Z_m(t)||\le K(1+2pl)||X(\zeta_i)|| {\rm e}^{-\alpha(t-\zeta_i)} < \frac{\epsilon}{2}, \, \zeta_i \le t,
\label{0909}
\end{eqnarray*}
The limit function $Z(t) = (X(t),Y(t))$ of the sequence is a solution of (\ref{treqn}) and satisfies
 \begin{eqnarray}
&& ||Z(t)|| \le K(1+2pl)||X(\zeta_i)|| {\rm e}^{-\alpha(t-\zeta_i)}< \frac{\epsilon}{2}, \, \zeta_i \le t,
\label{123909}
\end{eqnarray}
Hence, we can define a function $\tilde F(\zeta_i, X, d)$ such that $Y(\zeta_i) = \tilde F(\zeta_i, X(\zeta_i), d),$ which satisfies (\ref{F}).  Next, we can prove using (\ref{kryak}) and  (\ref{123**}) the existence of a pair
$(\bar X_0,\bar d)$ such that 
$$\bar X_0 = u_0  - G(\zeta_i,\bar d), v_0 - \bar d = \tilde F(\zeta_i, \bar X_0, \bar d),$$ 
$$||\bar X_0|| < \frac{\epsilon}{2K(1+ 2pl)},\,\, ||\bar d|| < 2 \delta.$$
Now,  transformation (\ref{tr}) and (\ref{123909}) imply that 
\begin{eqnarray}
&&||z(t, \zeta_i, z_0) - \mu((t, \zeta_i,\bar d)||\le K(1+2pl)||\bar X_0|| {\rm e}^{-\alpha(t-\zeta_i)}, \, \zeta_i \le t.
\label{hrup}
\end{eqnarray}

From (\ref{mutu}) and (\ref{hrup}) it follows that
\begin{eqnarray}
&&||z(t, \zeta_i, z_0)|| < \epsilon,\, \zeta_i \le t.
\label{hrupa}
\end{eqnarray}

Now, we can conclude in view of (\ref{hrupa}) that  the zero solution of 
(\ref{1}) is stable.
Assume that the zero solution of  (\ref{3*}) is  asymptotically stable, then (\ref{hrup}) implies that the zero solution of (\ref{1}) is also asymptotically stable.
Finally, it is  obvious that if  the zero solution of  (\ref{3*})  is unstable, then
the trivial solution of   (\ref{1}) is unstable as well.
The theorem is proved.

\vspace{0.5cm}
\begin{Ack}
The author wishes to express his sincere gratitude to the  referee for  the  helpful criticism and valuable suggestions, especially for the comment which encouraged the author to improve the  proof of the main result of the paper, Theorem \ref{9898}.

\end{Ack}

\end{document}